\documentclass{amsart}
\usepackage{amssymb}
\usepackage{amsthm}
\usepackage{graphics}
\newtheorem*{reftheorem}{Theorem}
\newtheorem{introtheorem}{Theorem}
\newtheorem{theorem}{Theorem}

\newtheorem{prop}{Proposition}

\numberwithin{equation}{section}
\title{Homotopically non-trivial maps with small k-dilation}
\author{Larry Guth}
\address{Department of Mathematics, Stanford, Stanford CA, 94305 USA}
\email{lguth@math.stanford.edu}

\begin{document}
\begin{abstract}
We construct homotopically non-trivial maps from $S^m$ to $S^n$
with arbitrarily small 3-dilation for certain pairs $(m,n)$.  The
simplest example is the case $m=4, n=3$, and there are other
pairs with arbitrarily large values of both m and n.  We show
that a homotopy class in $\pi_7(S^4)$ can be represented by maps
with arbitrarily small 4-dilation if and only if the class is
torsion.

\end{abstract}

\maketitle

The k-dilation of a map measures how much the map stretches
k-dimensional volumes.  If f is a $C^1$ map between Riemannian
manifolds, we say that the k-dilation of f is at most $D$ if f
maps each k-dimensional submanifold of the domain with volume V
to an image with volume at most $DV$.  We get the same k-dilation
whether we consider all submanifolds or whether we consider only
small disks, and so the k-dilation can also be defined in terms
of the first derivative df.  Recall that $\Lambda^k df$, the
k-fold exterior product of the derivative $df$, maps $\Lambda^k
TM$ to $\Lambda^k TN$.  If f is $C^1$, the k-dilation of f is
equal to the supremal value of the norm $|\Lambda^k df|$. 

In this paper, we examine to what extent a bound on the
k-dilation of a map controls the homotopy type of the map.  A
beautiful result of this kind was recently obtained by Tsui and
Wang in \cite{TW}.

\begin{reftheorem} (Tsui and Wang) Let f be a $C^1$ map from
$S^m$ to $S^n$, where $m \ge 2$.  If the 2-dilation of f is less
than 1, then f is nullhomotopic.
\end{reftheorem}

The main result of this paper shows that the situation is very
different for 3-dilation.

\begin{introtheorem} For each n, there are infinitely many m so that
the following holds: there are homotopically non-trivial maps
from $S^m$ to $S^n$ with arbitrarily small 3-dilation.
\end{introtheorem}

This result partly answers a question raised by Gromov in
\cite{G2} (page 231).  Gromov asked for which values of k, q, m,
and n is a map $f: S^m \rightarrow S^n$ with a sufficiently small
norm $|\Lambda^k df|_{L^q}$ necessarily null-homotopic.

We make the following definition.  A homotopy class $a$ in
$\pi_m(S^n)$ lies in $V_k \pi_m(S^n)$ if there are maps in the
homotopy class $a$ with arbitrarily small k-dilation.  We will
prove that $V_k \pi_m(S^n)$ is a subgroup of $\pi_m(S^n)$ and
that $V_k \pi_m(S^n) \subset V_{k+1} \pi_m(S^n)$.  Therefore,
$V_k \pi_m(S^n)$ defines a filtration of $\pi_m(S^n)$.  More
generally, we will define a filtration $V_k \pi_m(X)$ for any
space X.

Our methods give some partial information about the filtration
$V_k \pi_m(S^n)$.  The information is most interesting for the
filtration $V_k \pi_7(S^4)$.  Recall that $\pi_7(S^4)$ is
isomorphic to $\mathbb{Z} \oplus \mathbb{Z}_{12}$.  

\begin{introtheorem} The group $V_4 \pi_7(S^4)$ is the torsion
subgroup of $\pi_7(S^4)$.  It is a proper, non-zero subgroup.
\end{introtheorem}

The paper is organized as follows.  In the first section, we
prove the first theorem.  In the second section, we summarize
some lower bounds for k-dilation that appear in the literature. 
In the third section, we define $V_k \pi_m(X)$, and prove its
basic algebraic properties.  In the fourth section, we look at
some examples, including $V_k \pi_7(S^4)$.

This paper is based on a section of my thesis.  I would like to
thank my thesis advisor, Tom Mrowka, for his help and support.

\section{Homotopically non-trivial maps with small k-dilation}

This section contains the main result of the paper.  It gives a
construction for maps between spheres with very small k-dilation
which are homotopic to high order suspensions.

\begin{prop} Fix a homotopy class $a$ in $\pi_m(S^n)$ and then
consider its p-fold suspension $\Sigma^p a$ in
$\pi_{m+p}(S^{n+p})$, and let k be an integer greater than $n +
(n/m) p$.  Then there are maps from $S^{m+p}$ to $S^{n+p}$ in the
homotopy class $\Sigma^p a$ with arbitrarily small k-dilation.
\end{prop}

\proof Let $f_1$ be a map in the homotopy class $a$ from
$[0,1]^m$ to the unit n-sphere, taking the boundary of the domain
to the basepoint of $S^n$.  Let $f_2$ be a degree 1 map from
$[0,1]^p$ to the unit p-sphere, taking the boundary of the domain
to the basepoint of $S^p$.  We can assume both maps are $C^1$,
and we pick a number L which is bigger than the Lipshitz constant
of either map.

Inside of the unit (m+p)-sphere, we can quasi-isometrically embed
a rectangle R with dimensions $[0,\epsilon]^m \times
[0,\epsilon^{-m/p}]^p$.  (The quasi-isometric constant does not
depend on $\epsilon$.)  Now we construct a map F from R to $S^n
\times S^p$.  The map F is a direct product of a map $F_1$ from
$[0, \epsilon]^m$ to $S^n$ and a map $F_2$ from $[0,
\epsilon^{-m/p}]^p$ to $S^p$.  The map $F_1$ is just a dilation
from $[0, \epsilon]^m$ to the unit cube, composed with the map
$f_1$.  The map $F_2$ is just a dilation from $[0,
\epsilon^{-m/p}]^p$ to $[0,1]^p$, composed with the map $f_2$.  

When k is bigger than n, the k-dilation of F is less than $(L
\epsilon^{-1})^n (L \epsilon^{m/p})^{k-n}$.  Expanding this
expression gives $L^k \epsilon^{-n + (m/p)k - (m/p)n}$.  The
important part of the expression is the power of $\epsilon$,
which is equal to $(m/p) (k-n-(n/m)p)$.  We have assumed that k
is greater than $n + (n/m) p$, and so the exponent of $\epsilon$
is positive.  For $\epsilon$ sufficiently small, the k-dilation
of F is arbitrarily small.  

The map F takes the boundary of R to $S^n \vee S^p$.  We compose
F with a smash map, which is a degree 1 map from $S^n
\times S^p$ to $S^{n+p}$, taking $S^n \vee S^p$ to the base
point.  The result is a map from R to $S^{n+p}$ which takes the
boundary of R to the basepoint.  We can easily extend this map to
all of $S^{m+p}$ by mapping the complement of R to the basepoint
of $S^{n+p}$.  The resulting map is homotopic to $\Sigma^p(a)$,
and it has arbitrarily small k-dilation.
\endproof

We can apply our proposition to the non-trivial homotopy class in
$\pi_4(S^3)$, which is represented by the suspension of the Hopf
fibration.  In this case, we are considering a 1-fold suspension
of a map from $S^3$ to $S^2$.  Therefore $p=1$, $m=3$, and $n=2$. 
Since $3 > 2 + (2/3) 1$, this homotopy class can be realized by
maps from $S^4$ to $S^3$ with arbitrarily small 3-dilation.

We now turn to the proof of our main theorem, which gives
infinitely many examples of non-trivial homotopy classes that can
be realized with arbitrarily small 3-dilation.  The proof of this
theorem uses some deep results in algebraic topology to guarantee
that certain high-order suspensions are homotopically
non-trivial.  I would like to thank Haynes Miller for helping me
to find the relevant results in the topology literature.

\begin{theorem} For every $N \ge 2$, there are infinitely many M
so that there are homotopically non-trivial maps from $S^M$ to $S^N$
with arbitrarily small 3-dilation.
\end{theorem}

\proof When $i = 8j+1$, the homotopy group $\pi_i(SO)$ is equal
to $\mathbb{Z}_2$.  The stable J-homomorphism maps $\pi_i(SO)$ to
the $i^{th}$ stable stem of the homotopy groups of spheres.  The
image of the J-homomorphism was studied by Adams in \cite{A}. 
When $i=8j + 1$, he proved that the stable J-homomorphism is
injective.  Its image is a copy of $\mathbb{Z}_2$.  This image
contains a non-trivial element in $\pi_{i+n}(S^n)$, for large n. 
It turns out that this non-trivial element is the suspension of a
class in $\pi_{i+2}(S^2)$.  This statement is made clearly in the
introduction to the paper \cite{DM}, and the proof appears in the
older paper \cite{C}.  For each $i=8j+1$, let $a_i$ be a homotopy
class in $\pi_{i+2}(S^2)$ whose suspension is the non-trivial
element in the image $J(\pi_i(SO))$.  In particular, the p-fold
suspension $\Sigma^p a_i$ is non-trivial for every p and every i.

Now let N be any integer at least 2.  For each $i=8j+1$ greater
than $2N-6$, consider the class $\Sigma^{N-2} a_i$ in
$\pi_{i+N}(S^N)$.  We let $M = N+i$.  Each of these homotopy
classes is non-trivial.  In the language of Proposition 1, we
have $p = N-2$, $m=i+2$, and $n=2$.  The condition that $i$ is
greater than $2N - 6$ exactly guarantees that $3 > n + (n/m) p$. 
Therefore, each of these homotopy classes can be realized with
arbitrarily small 3-dilation.
\endproof

This proof constructs homotopically non-trivial maps from $S^M$
to $S^N$ with arbitrarily small 3-dilation for many pairs
$(M,N)$.  For example, when $N=3$, we can take $M=4, 12, 20, 28,
36,$ and so on; and when $N=4$, we can take $M = 13, 21, 29,
37,$ and so on.

\section{Lower bounds for k-dilation}

In this section, we survey several lower bounds for the
k-dilation of mappings.  We begin by recalling the theorem of
Tsui and Wang mentioned in the introduction.

\begin{reftheorem} (Tsui and Wang, \cite{TW}) Let f be a $C^1$
map from $S^m$ to $S^n$, where $m \ge 2$.  If the 2-dilation of f
is less than 1, then f is nullhomotopic.
\end{reftheorem}

The proof by Tsui and Wang uses the mean curvature flow to deform
the graph of the map f as a submanifold of $S^m \times S^n$. 
They prove that the mean curvature flow converges to the graph of
a constant function and that at each time t the flowed
submanifold is the graph of a map $f_t$.  Therefore, $f_t$
provides a homotopy from t to a constant map.

Earlier, Gromov had proved a slightly weaker theorem in the same
spirit.  He proved that for each m and n, there exists a number
$\epsilon(m,n) > 0$, so that any $C^1$ map from $S^m$ to $S^n$
with 2-dilation less than $\epsilon(m,n)$ is null-homotopic. 
Gromov's proof is completely different from the proof of Tsui and
Wang.  He approaches the map from $S^m$ to $S^n$ as a family of
maps from $S^2$ to $S^n$.  If the 2-dilation is less than
$\epsilon$, then each image of $S^2$ has area less than $4 \pi
\epsilon$.  Gromov uses the Riemann mapping theorem to change
coordinates on the domain so that each map in the family has
Dirichlet energy less than $C \epsilon$.  Then he uses the
borderline Sobolev inequality, which bounds the BMO norm of each
map in the family by $C \epsilon$.  Using the bound on the BMO
norm, he constructs a homotopy from the initial family of maps to
a family of constant maps.  This argument is sketched on page 179
of the long essay \cite{G2}.

Gromov's estimates are less sharp than those of Tsui and Wang,
but he is able to push them through with much weaker hypotheses. 
For example, Gromov proves that if f maps $S^m$ to $S^n$ with
$|\Lambda^2 df|_{L^{m-1+\delta}} < \epsilon(m,n,\delta)$ then f
is null-homotopic.  Gromov's method also applies to more
manifolds.  He proves that for any two compact simply connected
Riemannian manifolds M and N there is a constant $\epsilon(M,N)$
so that any map $f:M \rightarrow N$ with 2-dilation less than
$\epsilon(M,N)$ is null-homotopic.  (Both these results are
special cases of the Corollary on page 230 of \cite{G2}.)

Lower bounds for the k-dilation with k greater than 2 are much
rarer.  We begin with a very elementary example.  If a $C^1$ map
f from $S^n$ to $S^n$ is homotopically non-trivial, then it has
n-dilation at least 1.  This result follows because a map with
n-dilation equal to $1-\epsilon$ has an image with volume at most
$(1-\epsilon) Vol(S^n)$.  Such a map is not surjective, and so it
is null-homotopic.

Gromov also proved lower bounds for the k-dilation of maps with
non-zero Hopf invariant.  The following argument appears in
\cite{G1} on pages 358-359.

\begin{reftheorem} (Gromov) Let f be a map from $S^{4n-1}$ to
$S^{2n}$ with $2n$-dilation D.  Then the norm of the Hopf
invariant of f is bounded by $C(n) D^2$.  Since the Hopf
invariant is an integer, any map with non-zero Hopf invariant has
$2n$-dilation at least $C(n)^{-1/2}$.
\end{reftheorem}

\proof Let $\omega$ be a 2n-form on $S^{2n}$ with $\int \omega =
1$.  The pullback $f^*(\omega)$ is a closed 2n-form on
$S^{4n-1}$.  Since $H^{2n}(S^{4n-1}) = 0$, this form is exact. 
We let $P f^*(\omega)$ denote any primitive of $f^*(\omega)$. 
Then the Hopf invariant of f is equal to $\int_{S^{4n-1}} P
f^*(\omega) \wedge f^*(\omega)$.

During this proof, we let $C$ denote a constant depending on n
that may change from line to line.  We can assume that $|\omega|
< C$ at every point of $S^{2n}$.  The norm of $f^*(\omega)$ is
bounded by $C D$ pointwise.  Therefore, the $L^2$ norm of
$f^*(\omega)$ is bounded by $C D$.  Using Hodge theory, we can
choose $P f^*(\omega)$ to be perpendicular to all of the exact
(2n-1)-forms.  For this choice, the $L^2$ norm of $P f^*(\omega)$
is bounded by $\lambda^{-1/2} C D$, where $\lambda$ is the
smallest eigenvalue of the Laplacian on exact (2n)-forms.  The
eigenvalue $\lambda$ is greater than zero and depends only on n. 
Finally, the norm of the Hopf invariant is bounded by
$|f^*(\omega)|_{L^2} |P f^*(\omega)|_{L^2}$, which is bounded by
$C(n) D^2$. \endproof

The method above can be generalized to many non-torsion homotopy
classes.  For more information, see Gromov's discussion in
\cite{G2}, pages 220-223.

The methods outlined here leave many open questions.  If $a$ is a
torsion homotopy class in some homotopy group $\pi_m(S^n)$, I
don't know any way to lower bound the 3-dilation of maps in the
homotopy class $a$.  For example, the 100-fold suspension of the
Hopf map is a non-trivial element in $\pi_{103}(S^{102})$. 
According to Proposition 1, it can be realized by maps
with arbitrarily small 69-dilation.  It would be interesting to
know whether it can be realized by maps with arbitrarily small
3-dilation.

\section{A filtration on homotopy groups}

In this section, we define a filtration $V_k \pi_m(X)$ which
measures which homotopy classes in $\pi_m(X)$ can be realized
with arbitrarily small k-dilation.  We derive some easy formal
properties of this filtration.  In the next section, we will
calculate it in some examples and show that it can be
non-trivial.

As a first step, we define $V_k \pi_m(W)$ when W is a finite
simplicial complex.  We put a piecewise Riemannian metric on W by
equipping each simplex with the metric of an equilateral simplex
in Euclidean space with edges of length 1.  We then define $V_k
\pi_m(W)$ to be the set of homotopy classes in $\pi_m(W)$ that
can be realized by maps $S^m \rightarrow W$ with arbitrarily
small k-dilation.

The set $V_k \pi_m(W)$ is in fact a subgroup of $\pi_m(W)$.  If
$a$ lies in $V_k \pi_m(W)$, then let $f_i$ be a sequence of maps
from $S^m$ to $W$ in the homotopy class $a$ with k-dilation
tending to zero.  Let $I$ be a reflection, mapping $S^m$ to
itself with degree -1, and taking the basepoint of $S^m$ to
itself.  Then the maps $f_i \circ I$ have k-dilations tending to
zero and lie in the homotopy class $-a$.  Therefore $-a$ lies in
$V_k \pi_m(W)$.  Next, suppose that $a$ and $b$ lie in $V_k
\pi_m(W)$.  Again, let $f_i$ be a sequence of (pointed) maps in
the class $a$ with k-dilation tending to zero, and let $g_i$ be a
sequence of (pointed) maps in the homotopy class $b$ with
k-dilation tending to zero.  Let $I$ be a map from $S^m$ to $S^m
\vee S^m$ with degree (1,1).  Let $h_i$ be the map from $S^m \vee
S^m$ to W whose restriction to the first copy of $S^m$ is equal
to $f_i$ and whose restriction to the second copy of $S^m$ is
equal to $g_i$.  Then the sequence $h_i \circ I$ has k-dilation
tending to zero.  Each map in the sequence lies in the homotopy
class $a+b$.  So $a+b$ lies in $V_k \pi_m(W)$.

The sets $V_k \pi_m(W)$ are nested, with $V_k \pi_m(W) \subset
V_{k+1} \pi_m(W)$.  To see this, we express the k-dilation of a
map in terms of the singular values of its derivative.  Let f be
a piecewise $C^1$ map from $S^m$ to W.  Suppose that the
singular values of $df$ at a point are equal to $0 \le s_1 \le ...
\le s_m$.  Then the norm $|\Lambda^k df|$ at that point is equal
to $s_{m-k+1} ... s_m$.  Therefore, $|\Lambda^{k+1} df| \le
|\Lambda^k df|^{\frac{k+1}{k}}$.  If a map f has k-dilation D,
then the (k+1)-dilation of f is at most $D^{\frac{k+1}{k}}$.  In
particular, if $f_i$ is a sequence of maps with k-dilation
tending to zero, then the (k+1)-dilation of $f_i$ also tends to
zero.  Therefore, $V_k \pi_m(W) \subset V_{k+1} \pi_m(W)$.  Any
map from $S^m$ has $(m+1)$-dilation zero, and so $V_{m+1}
\pi_m(W)$ is always the whole homotopy group $\pi_m(W)$.

The last paragraph shows that the sequence of groups $V_k
\pi_m(W)$ is a filtration of $\pi_m(W)$.

$ 0 = V_1 \pi_m(W) \subset V_2 \pi_m(W) \subset ... \subset V_m
\pi_m(W) \subset V_{m+1} \pi_m(W) = \pi_m(W).$

The filtration $V_k$ behaves naturally under mappings in the
following sense.  If $F: W \rightarrow V$ is a continuous pointed
mapping between finite simplicial complexes, then $F_*: \pi_m(W)
\rightarrow \pi_m(V)$ takes $V_k \pi_m (W)$ into $V_k \pi_m(V)$. 
To prove this, first approximate F by a PL map with some finite
Lipshitz constant L.  Let $a$ be a class in $V_k \pi_m(W)$,
realized by mappings $f_i: S^m \rightarrow W$ with k-dilation
tending to zero.  The map $F \circ f_i$ from $S^m$ to V has
k-dilation less than $L^k$ times the k-dilation of $f_i$, so the
the sequence $F \circ f_i$ has k-dilation tending to zero.  Each
map $F \circ f_i$ lies in the homotopy class $F_*(a)$. 
Therefore, $F_*(a)$ lies in $V_k \pi_m(V)$.

Now we define $V_k \pi_m(X)$ for an arbitrary space X.  A
homotopy class $a$ in $\pi_m(X)$ belongs to $V_k \pi_m(X)$ if
there is a finite simplicial complex W and a map $F: W
\rightarrow X$ so that $a$ lies in the image $F_*(V_k \pi_m(W))$. 
In case X is a finite simplicial complex, this definition agrees
with our first definition because of the mapping property of
$V_k$ proved in the last paragraph.  As above, $V_k \pi_m(X)$
defines a filtration of $\pi_m(X)$.  Also, if $f: X \rightarrow
Y$ is any continuous map of spaces, then $f_*$ maps $V_k
\pi_m(X)$ into $V_k \pi_m(Y)$.  From this mapping property, it
follows that the filtration $V_k \pi_m(X)$ is a homotopy
invariant.

We can rephrase the main theorems of this paper in the language
of the filtration $V_k \pi_m(S^n)$.  The theorem of Tsui and Wang
implies that $V_2 \pi_m(S^n) = 0$ for $m \ge 2$.  Theorem 1 says
that for each $n \ge 2$, the group $V_3 \pi_m(S^n)$ is non-zero for
infinitely many m.  Theorem 2 says that $V_4 \pi_7(S^4)$ is
exactly the torsion subgroup of $\pi_7(S^4)$, which is isomorphic
to $\mathbb{Z}_{12}$.

\section{Examples of the $V_k$ Filtration}

In this section, we use the tools from sections 1 and 2 to
compute the filtration $V_k \pi_m(S^n)$ for a few small values of
m and n.  For other values of m and n, we obtain some partial
information about the filtration.  The most interesting example
is a computation of $V_4 \pi_7 (S^4)$.  We use the lists of
homotopy groups of spheres and the suspension maps between them
given in \cite{T} on pages 39-42.  According to the the theorem of
Tsui and Wang, $V_2 \pi_m(S^n)$ is zero in all the examples
below, so we will only discuss $V_k$ for $k \ge 3$.

{\bf The filtration of $\pi_m(S^2)$}

Any element of $\pi_m(S^2)$ lies in $V_3 \pi_m(S^2)$, because any
smooth map from $S^m$ to $S^2$ has 3-dilation equal to zero.

{\bf The filtration of $\pi_4(S^3)$}

The group $\pi_4(S^3)$ is isomorphic to $\mathbb{Z}_2$, and the
non-trivial element is given by the suspension of the Hopf
fibration from $S^3$ to $S^2$.  Applying Proposition 1, this
non-trivial element lies in $V_3 \pi_4(S^3)$.

{\bf The filtration of $\pi_5(S^3)$}

The group $\pi_5(S^3)$ is isomorphic to $\mathbb{Z}_2$, and the
non-trivial element is the suspension of a class in $\pi_4(S^2)$. 
Applying Proposition 1, this non-trivial element lies in $V_3
\pi_5(S^3)$.

{\bf The filtration of $\pi_6(S^3)$}

The group $\pi_6(S^3)$ is isomorphic to $\mathbb{Z}_{12}$.  One
non-trival element is the suspension of a class in $\pi_5(S^2)$. 
Applying Proposition 1, this non-trivial element lies in $V_3
\pi_6(S^3)$.  I don't know whether the other non-trivial elements
lie in $V_3 \pi_6(S^3)$.

{\bf The filtration of $\pi_5(S^4)$}

The group $\pi_5(S^4)$ is isomorphic to $\mathbb{Z}_2$, and the
non-trivial element is given by the double suspension of the Hopf
fibration from $S^3$ to $S^2$.  Applying Proposition 1, this
non-trivial element lies in $V_4 \pi_5(S^4)$.  I don't know
whether it lies in $V_3 \pi_5(S^4)$.

{\bf The filtration of $\pi_6(S^4)$}

The group $\pi_6(S^4)$ is isomorphic to $\mathbb{Z}_2$, and the
non-trivial element is the double suspension of a class in
$\pi_4(S^2)$.  Applying Proposition 1, this non-trivial element
lies in $V_4 \pi_6(S^4)$.  I don't know whether it lies in $V_3
\pi_6(S^4)$.

{\bf The filtration of $\pi_7(S^4)$}

The group $\pi_7(S^4)$ is isomorphic to $\mathbb{Z} \oplus
\mathbb{Z}_{12}$.  The Hopf invariant gives a map $H: \pi_7(S^4)
\rightarrow \mathbb{Z}$.  In section 2, following Gromov, we
proved that $V_4 \pi_7(S^4)$ lies in the kernel of H. 
The kernel of H is isomorphic to $\mathbb{Z}_{12}$.  The
suspension map is an isomophism from $\pi_6(S^3)$ onto the kernel
of H, as follows from the long exact sequence of the Hopf
fibration $S^3 \rightarrow S^7 \rightarrow S^4$.  (This fact is
stated in \cite{T} on page 2.)  Therefore, every element in the
kernel of H is a suspension.  Applying Proposition 1, we see that
the kernel of H lies in $V_4 \pi_7(S^4)$.  Therefore, $V_4
\pi_7(S^4)$ is exactly the kernel of H.  This proves the
second theorem in the introduction.

In addition, one non-trivial element in $\pi_7(S^4)$ is a double
suspension of an element in $\pi_5(S^2)$.  Applying Proposition 1
to this element, we see that it lies in $V_3 \pi_7(S^4)$.  I
don't know whether the other torsion elements of $\pi_7(S^4)$ lie
in $V_3 \pi_7(S^4)$.

\end{document}